\documentclass{article}
\usepackage{amssymb,epsf,latexsym,a4}
\usepackage[english]{babel}
\newcommand{\qed}{\unskip\nobreak\hfill\mbox{ $\Box$}\par}
\def\bbbr{{\mathbb R}}
\def\bbbc{{\mathbb C}}

\def\bbbz{{\mathbb Z}}

\def\mod#1{({\rm mod\ }#1)}
\def\v#1{{\bf#1}}

\def\D{\partial}

\def\dim{{\rm dim}}
\def\supp{{\rm supp}}

\def\is{\equiv}
\def\mod#1{({\rm mod}\ #1)}
\newtheorem{theorem}[subsection]{Theorem}
\newtheorem{definition}[subsection]{Definition}
\newtheorem{lemma}[subsection]{Lemma}

\newtheorem{corollary}[subsection]{Corollary}

\newtheorem{proposition}[subsection]{Proposition}

\title{Irreducibility of A-hypergeometric systems}
\author{F.Beukers}
\date{September 1, 2010}

\begin{document}
\maketitle

\abstract{We give an elementary proof of the Gel'fand-Kapranov-Zelevinsky
theorem that non-resonant A-hypergeometric
systems are irreducible. We also provide a proof of a converse statement}

\section{Introduction}
Let $A\subset \bbbz^r$ be a finite set such that
\begin{enumerate}
\item The $\bbbz$-span of $A$ is $\bbbz^r$.
\item There exists a linear form $h$ such that $h(\v a)=1$ for
all $\v a\in A$.
\end{enumerate}
Let $\alpha=(\alpha_1,\ldots,\alpha_r)\in\bbbc^r$. At the end of the
1980's Gel'fand, Kapranov and Zelevinsky \cite{GKZ1}, \cite{GKZ2}, \cite{GKZ3}
developed a theory of hypergeometric
functions and equations which uses $A$ and $\alpha$ as starting data. It turns
out that the resulting equations contain the classical cases of Appell, Horn,
Lauricella and Aomoto hypergeometric functions.

Denote $A=\{\v a_1,\ldots,\v a_N\}$ (with $N>r$). Writing the vectors $\v a_i$ in
column form we get the so-called A-matrix
$$A=\pmatrix{
a_{11} & a_{12} & \cdots & a_{1N}\cr
a_{21} & a_{22} & \cdots & a_{2N}\cr
\vdots &        &        & \vdots\cr
a_{r1} & a_{r2} & \cdots & a_{rN}\cr
}
$$
For $i=1,2,\ldots,r$ consider the first order differential operators
$$Z_i=a_{i1}v_1\D_1+a_{i2}v_2\D_2+\cdots+a_{iN}v_N\D_N$$
where $\D_j={\D\over \D v_j}$ for all $j$.

Let
$$L=\{(l_1,\ldots,l_N)\in\bbbz^N|\ l_1\v a_1+l_2\v a_2+\cdots+l_N\v a_N=\v 0\}$$
be the lattice of integer relations between the elements of $A$. For every $\v l\in L$
we define the so-called box-operator
$$\Box_{\v l}=\prod_{l_i>0}\D_i^{l_i}-\prod_{l_i<0}\D_i^{-l_i}$$
The system of differential equations
\begin{eqnarray*}
(Z_i-\alpha_i)\Phi&=&0\quad (i=1,\ldots,r)\\
\Box_{\v l}\Phi&=&0\quad \v l\in L
\end{eqnarray*}
is known as the system of {\it A-hypergeometric differential equations} and we denote it
by $H_A(\alpha)$. We like to remark that
independently, and at around the same time, B.Dwork arrived at a similar setup for generalised
hypergeometric functions. The system of A-hypergeometric equations is implicit in his book \cite{dwork}.

Let $K=\bbbc(v_1,\ldots,v_N)$ and let ${\cal H}_A(\alpha)$ be the left ideal in $K[\D_1,\ldots,\D_N]$
generated by the operators from $H_A(\alpha)$. The quotient $K[\D_1,\ldots,\D_N]/{\cal H}_A(\alpha)$
is a $K$-module. Its $K$-rank is called the rank of the system $H_A(\alpha)$. Furthermore,
the system is called {\it non-resonant} if the set $\alpha+\bbbz^r$ has empty intersection with the
boundary of $C(A)$. The system is called {\it resonant} if the intersection is non-empty.

In \cite{GKZ3}, (corrected in \cite{GKZ5}) and \cite[Corollary 5.20]{adolphson} the following theorem is shown.

\begin{theorem}[GKZ, Adolphson] Suppose either one of the following conditions holds,
\begin{enumerate}
\item the toric ideal $I_A$ in $\bbbc[\D_1,\ldots,\D_N]$ generated by the box
operators has the Cohen-Macaulay property.
\item The system $H_A(\alpha)$ is non-resonant.
\end{enumerate}
Then the rank of $H_A(\alpha)$ is finite and equals
the volume of the convex hull $Q(A)$ of the points of $A$. The volume is normalized so that a minimal
$(r-1)$-simplex with integer vertices in $h(\v x)=1$ has volume $1$.
\end{theorem}

Let $p$ be a generic point in $(\bbbc^*)^N$ (the space with coordinates $v_1,\ldots,v_N$). Then it is
known that the dimension of the $\bbbc$-vector space of local power series solutions around $p$ of $H_A(\alpha)$ equals
the rank of $H_A(\alpha)$.

The $K$-module $K[\D_1,\ldots,\D_N]/{\cal H}_A(\alpha)$ has a natural left action by the operators $\D_i$, so it is a
D-module. We shall say that the system $H_A(\alpha)$ is {\it irreducible} if this D-module has no
submodules beside $0$ and the module itself. We call it {\it reducible} otherwise.
Gel'fand, Kapranov and Zelevinsky proved in \cite[Thm 2.11]{GKZ4} the following beautiful theorem.

\begin{theorem}[GKZ, 1990]\label{nonresonant2irreducible}
Suppose the system $H_A(\alpha)$ is non-resonant.
Then $H_A(\alpha)$ is irreducible.
\end{theorem}

The proof uses the theory of perverse sheaves and is hard to follow for someone without this background.
It is the purpose of the present paper to give a more elementary proof of this theorem.
This is done in Section \ref{proofn2i}.
In addition we prove a converse statement, namely the following.

\begin{theorem}\label{resonant2reducible} Suppose that the toric ideal $I_A$ has the Cohen-Macaulay property
and suppose that the convex hull $Q(A)$
is not a pyramid. If the system $H_A(\alpha)$ is resonant, then it is reducible.
\end{theorem}

As far as we could see the latter theorem is not stated as such in the papers of Gel'fand, Kapranov and Zelevinsky
or any other papers.
The condition that $Q(A)$ is not a pyramid means that we like to avoid the situation where $A$ contains $N-1$
points in an $r-2$-dimensional affine plane and only one outside of it. It is not hard to see that $Q(A)$ is a
pyramid if and only if for every index $i\in\{1,\ldots,N\}$ there exists $\v l\in L$ such that $l_i\ne0$.
Suppose $Q(A)$ is a pyramid with top $\v a_1$. Then one easily sees that the box-operators do not contain $\D_1$.
Hence there exists $\beta\in\bbbr$ such that the solutions
of $H_A(\alpha)$ have the form $v_1^{\beta}F(v_2,\ldots,v_N)$. In case $\beta=0$, so all solutions independent
of $v_1$, the vector of parameters lies in the bottom of the pyramid, which is the affine space spanned
by $\v a_2,\ldots,\v a_N$.


\section{Contiguity}
Consider the system $H_A(\alpha)$,
$$\Box_{\v l}\Phi=0,\ \v l\in L,\qquad Z_j\Phi=\alpha_j\Phi,\ j=1,\ldots,r.$$
Apply the operator $\D_i$ from the left. We obtain,
$$\Box_{\v l}\D_i\Phi=0,\ \v l\in L,
\qquad Z_j\D_i\Phi=-a_{ji}\D_i\Phi,\ j=1,\ldots,r.$$
In other words, $F\mapsto \D_iF$ maps the solution space of $H_A(\alpha)$
to the solution space of $H_A(\alpha-\v a_i)$.

We can phrase this alternatively in terms of D-modules. Denote by ${\cal H}_A(\alpha)$
the left ideal in $K[\D]$ generated by the hypergeometric operators
$\Box_{\v l}$ and $Z_j$. Then the map $P\mapsto P\D_i$ gives a D-module homomorphism
$K[\D]/{\cal H}_A(\alpha-\v a_i)\to K[\D]
/{\cal H}_A(\alpha)$. We are interested in the cases when this is a D-module isomorphism or,
equivalently, whether $F\mapsto\D_iF$ gives an isomorphism of solution spaces.

The following Theorem was first proven by B.Dwork in his book \cite[Thm 6.9.1]{dwork}. Another
proof was given in \cite[Lemma 7.10]{dicmatmil}.
We present an adaptation of Dwork's ideas into a language which is quite different from Dwork's.

\begin{theorem}[Dwork]\label{dwork}
Suppose $H_A(\alpha)$ is non-resonant. Then the map $F\mapsto\D_iF$ yields an isomorphism
between the solution spaces of $H_A(\alpha)$ and $H_A(\alpha-\v a_i)$.
\end{theorem}

For the proof we need an extra Lemma and some notation.
Suppose the positive cone $C(A)$ is given by a finite set ${\cal F}$ of linear inequalities $l(\v x)\ge0,
\ l\in{\cal F}$. Assume moreover that the linear forms $l$ are integral valued on $\bbbz^r$ and
normalise them so that the greatest common divisor of all values is $1$.

Consider the integral points in $C(A)$. It is not necessarily true that every point in
$C(A)\cap\bbbz^r$ is a linear combination of the $\v a_i$ with non-negative integer
coefficients. However, we do have the following Lemma.

\begin{lemma}\label{saturation}
There exists a point $p\in C(A)\cap\bbbz^r$ such that $(\v p+C(A))\cap\bbbz^r\subset\bbbz_{\ge0}A$
where $\bbbz_{\ge0}A$ is the span of $A$ with non-negative integer coefficients.
\end{lemma}

{\bf Proof}
It is clear that there exists a positive integer $\delta$ such that for any point $(\lambda_1,\ldots,\lambda_N)\in
L\otimes\bbbr$ there exists $(m_1,\ldots,m_n)\in L$ such that $|m_i-\lambda_i|\le\delta$.
Let us take $\v p=\delta(\v a_1+\cdots+\v a_N)$.

Suppose we are given a point $\v n\in (\v p+C(A))\cap\bbbz^r$. Then there exist $\lambda_i\in\bbbr_{\ge\delta}$
and integers $n_1\ldots,n_N$ such that $\v n=\lambda_1\v a_1+\cdots+\lambda_N\v a_N
=n_1\v a_1+\cdots+n_N\v a_N$. The point $(\lambda_1-n_1,\ldots,\lambda_N-n_N)$ lies in
$L\otimes\bbbr$. Hence there exists $(m_1,\ldots,m_N)\in L$ such that $|\lambda_i-n_i-m_i|\le
\delta$ for $i=1,\ldots,N$. Since $\lambda_i\ge\delta$ for every $i$ we find that
$n_i+m_i\ge0$. Hence $\v n=n_1\v a_1+\cdots+n_N\v a_N=(n_1+m_1)\v a_1+\cdots+(n_N+m_N)\v a_N$,
hence $\v n\in\bbbz_{\ge0}A$.
\qed
\medskip

{\bf Proof} of Thm \ref{dwork}.
We will construct an operator $P\in K[\D]$
such that $P\D_i\is 1\mod{{\cal H}_A(\alpha)}$. In particular, $F\mapsto P(F)$
would be the inverse of $\D_i$, which establishes the isomorphism.

For any $l\in {\cal F}$ and any differential operator $\D^{\v u}=\D_1^{u_1}\cdots\D_N^{u_N}$
we define the valuation
$val_l(\D^{\v u})=\sum_{j=1}^Nu_jl(\v a_j)$. More generally, for any differential
operator $P\in K[\D]$ we define $val_l(P)$ to be the minimal valuation
of all terms in $P$.

Let $\v p$ be as in Lemma \ref{saturation}.
Suppose $val_l(\D^{\v u})\le val_l(\D^{\v w})+l(\v p)$ for every $l\in{\cal F}$. Hence
$\sum_{j=1}^N l((w_j-u_j)\v a_j)\ge l(\v p)$ for all $l\in{\cal F}$. So, according to
Lemma \ref{saturation} $\sum_{j=1}^N(w_j-u_j)\v a_j$ is a lattice
point in $\bbbz_{\ge0}A$. Hence there exist non-negative integers $w'_j$
such that $\sum_{j=1}^Nw'_j\v a_j=\sum_{j=1}^N(w_j-u_j)\v a_j$. Hence $\D^{\v w}$ is equivalent
modulo the box operator $\Box_{\v w-\v w'-\v u}$ with $\D^{\v w'}\D^{\v u}$.

Let $l\in{\cal F}$ be given. We show that modulo the ideal ${\cal H}_A(\alpha)$, the operator
$\D^{\v u}$ is equivalent to an operator $P$ such that $val_l(P)>val_l(\D^{\v u})$ and
$v_{l'}(P)\ge v_{l'}(\D^{\v u})$ for all $l'\in{\cal F},\ l'\ne l$. Let
$Z_l=-l(\alpha)+\sum_{j=1}^Nl(\v a_j)v_j\D_j$. Notice that $Z_l\in{\cal H}_A(\alpha)$ and $\D^{\v u}Z_l=
Z_l\D^{\v u}+l(\v u)\D^{\v u}$. Hence,
$$\sum_{j=1}^Nl(\v a_j)v_j\D_j\D^{\v u}\is l(\alpha-\v u)\D^{\v u}
\mod{{\cal H}_A(\alpha)}.$$
For each term on the left we have $l(\v a_j)\ne0\Rightarrow val_l(\D_j\D^{\v u})
>val_l(\D^{\v u})$.
Since, by non-resonance, $l(\alpha-\v u)\ne0$ our assertion is proven. Choose $k_l\in\bbbz_{\ge0}$
for every $l\in{\cal F}$. By repeated application
of our principle we see that any monomial $\D^{\v u}$ is equivalent modulo ${\cal H}_A(\alpha)$
to an operator $P$ with $val_l(P)\ge k_l+val_l(\D^{\v u})$ for all $l\in{\cal F}$.

In particular, there exists an operator $P$, equivalent to $1$ and $val_l(P)\ge val_l(\D_i)+l(\v p)$ for
every $l\in{\cal F}$. Then, $P$ is equivalent to an operator $P'\D_i$. Summarising,
$1\is P'\D_i\mod{{\cal H}_A(\alpha)}$. So $F\mapsto\D_iF$ is injective on the
solution space of $H_A(\alpha)$.

\qed
\medskip

There is another instance when $F\mapsto\D_iF$ is an isomorphism of solution spaces.
\begin{theorem}\label{irred2isom}
Suppose that the toric ideal $I_A$ has the Cohen-Macaulay property, that $Q(A)$ is not a pyramid
and that $H_A(\alpha)$ is an irreducible system.
Then $F\mapsto\D_iF$ gives an
isomorphism of solution spaces of $H_A(\alpha)$ and $H_A(\alpha-\v a_i)$.
\end{theorem}

{\bf Proof}.
Since $H_A(\alpha)$ is irreducible, the kernel of $F\mapsto\D_iF$ is either trivial
or the entire solution space. In the first case we are done, the map is injective
and the solution spaces have the same dimension (because $I_A$ has the Cohen-Macaulay property).

Now suppose we are in the second case, when $\D_iF\equiv 0$
for every solution $F$ of $H_A(\alpha)$. This is equivalent
to the statement $\D_i\in{\cal H}_A(\alpha)$. Let us write
$$\D_i=\sum_{\lambda}A_{\lambda}\Box_{\lambda}+\sum_{j=1}^rB_j(Z_j-\alpha_j).$$
The summation over the $\lambda\in L$ is supposed to be a finite summation. Let
us assume that we have chosen the $A_{\lambda}$ and $B_i$ such that the maximum
of the orders of the $B_i$ is minimal. Call this minimum $m$. We assert that $m=0$.
Suppose $m>0$.

We now work over the polynomial ring $R=\bbbc(\v v)[X_1,\ldots,X_N]$.
For any differential operator $P$ we write $P(\v X)$ for the polynomial we get after
we replace $\partial_j$ by $X_j$ for all $j$ in $P$. Write $I_A$ for the ideal in $R$
generated by the $\Box_{\v l}(\v X)$. Since the quotient ring
$R/I_A$ is a Cohen-Macaulay ring, the linear forms $Z_i(\v X)$ form a regular
sequence. In particular this means that if $P_1Z_1(\v X)+\cdots+P_rZ_r(\v X)=0$ in
$R/I_A$, then there exist polynomials $\eta_{ij}$ with $\eta_{ij}=-\eta_{ji}$
such that $P_i=\sum_{j=1}^r\eta_{ij}Z_j(\v X)$ for $i=1,\ldots,r$.

Let us return to the $A_{\lambda}$ and $B_j$ above. Note that $(A_{\lambda}\Box_{\lambda})(\v X)
=A_{\lambda}(\v X)\Box_{\lambda}(\v X)$ since the box-operators have constant coefficients.
Denote the order $m$ part of each $B_j$ by $B_j^{(m)}$. Then the $m+1$-st degree part of
$\sum_j (B_j(Z_j-\alpha_i))(\v X)$ reads $\sum_j B_j^{(m)}(\v X)Z_j(\v X)$. Since $m+1>1$
this degree $m+1$ part is zero in $R/I_A$. Hence there exist polynomials $\eta_{jk}$ with
$\eta_{jk}=-\eta_{kj}$ such that $B_j^{(m)}(\v X)=\sum_{k=1}^r\eta_{jk}Z_k(\v X)$ in $R/I_A$.
Denote bij $E_{jk}$ the differential operator which we get after we replace the variables $X_b$ in
$\eta_{jk}$ bij their counterparts $\D_b$. Define $\tilde{B}_j=B_j-\sum_{k=1}^rE_{jk}(Z_k-\alpha_k)$
and note that $\tilde{B}_j$ has order $<m$. Moreover,
$$\sum_{j=1}^rB_j(Z_j-\alpha_j)=\sum_{j=1}^r\tilde{B}_j(Z_j-\alpha_j)+\sum_{j,k=1}^rE_{jk}(Z_j-\alpha_j)
(Z_k-\alpha_k).$$
The last sum, by virtue of the antisymmetry of the $E_{jk}$ and the fact that $Z_j-\alpha_j$ and
$Z_k-\alpha_k$ commute for all $j,k$, is equal to zero in $R/I_A$.
Hence
$$\D_i\is \sum_{j=1}^r\tilde{B}_j(Z_j-\alpha_j)\mod{I_A}$$
where the $\tilde{B}_i$ have order $<m$. This contradicts the minimality of $m$.
Therefore we conclude that $m=0$. In other words there exist $b_i\in\bbbc(\v v)$ such
that $\D_i\equiv\sum_{j=1}^rb_j(Z_j-\alpha_j)\mod{I_A}$. Since the box-operators all
have order $\ge2$ this relation holds exact. It follows that there exist $\beta_j\in\bbbc$
such that $v_i\D_i=\sum_{j=1}^r\beta_j(Z_j-\alpha_j)$. In other words there exists a
linear form $m$ on $\bbbr^r$ such that $m(\v a_j)=0$ for all $j\ne i$ and $m(\v a_i)=1$.
But this implies that $Q(A)$ is a pyramid with $\v a_i$ as a top.
\qed

\section{Resonant systems}
In this section we prove Theorem \ref{resonant2reducible}.
Suppose that $H_A(\alpha)$ is resonant and irreducible. Then, by Theorem \ref{irred2isom}
for any $i$ the map $F\mapsto\D_iF$ is an isomorphism of solution spaces of
$H_A(\alpha)$ and $H_A(\alpha-\v a_i)$.
So we see that $H_A(\beta)$ is irreducible for any $\beta\in\bbbr^r$ with
$\beta\is\alpha\mod{\bbbz^r}$. Since the system is resonant there exists such a $\beta$ in a
face $F$ of $C(A)$. Suppose $A\cap F=\{\v a_1,\ldots,\v a_t\}$. We assert that there exist
non-trivial solutions of the form $f=f(v_1,\ldots,v_t)$. Suppose that $s={\rm rank}(\v a_1,
\ldots,\v a_t)$. By an $SL(r,\bbbz)$ change of coordinates we can
see to it that $F$ is given by $x_{s+1}=\cdots=x_r=0$. Then the coordinate $a_{rj}$ of $\v a_j$ is zero
for $i=s+1,\ldots,r$ and $j=1,\ldots,t$. Also, $\beta_{s+1}=\cdots=\beta_r=0$.
A solution $f=f(v_1,\ldots,v_t)$ satisfies the homogeneity equations
$$\left(-\beta_i+\sum_{j=1}^ta_{ij}v_j\D_j\right)f=0,\ i=1,\ldots,s.$$
Notice that the homogeneity equation with $i=s+1,\ldots,r$ are trivial.

Consider the box-operator $\Box_{\lambda}$ with $\lambda\in L$.
Write $\lambda=(\lambda_1,\ldots,\lambda_N)$. The positive support is the set
of indices $i$ where $\lambda_i>0$, the negative support is the set of
indices $i$ where $\lambda_i<0$.

Suppose the positive support is contained in $1,2,\ldots,t$. Then
$\sum_{\lambda_i>0}\lambda_i\v a_i$ is in ${\cal F}$. Hence
$-\sum_{\lambda_i<0}\lambda_i\v a_i$ is also in $F$. Since $F$
is a face, all non-zero terms of the latter have index $\le t$. So the negative
support is also in $1,2,\ldots,t$. Hence
$${\rm negative\ support}\subset\{1,\ldots,t\}
\iff{\rm positive\ support}\subset\{1,\ldots,t\}.$$
If the positive and negative support of $\lambda$ contain indices $>t$ then $f(v_1,\ldots,v_t)$
satisfies $\Box_{\lambda}f=0$ trivially.

Define a new set $\tilde{A}=\{\tilde{\v a}_1,\ldots,\tilde{\v a}_t\}\subset\bbbz^s$
where $\tilde{\v a}_j$ is the projection of $\v a_j$ on its first $s$ coordinates.
Define a new parameter $\tilde{\beta}$ similarly.
The solutions of the form $f(v_1,\ldots,v_t)$ of the original GKZ-system
satisfy the new GKZ-system corresponding
to $H_{\tilde{A}}(\tilde{\beta})$. They all satisfy the additional equations $\partial_iF=0$
for $i>t$, so they form a proper subspace of the solution space of $H_A(\alpha)$.
Hence the system is reducible, contradicting
our initial assumption of irreducibility.
\qed

\section{Series solutions}
Just as in the classical literature we like to be able to display explicit series solutions for the
A-hypergeometric system. In GKZ-theory one chooses $\gamma=(\gamma_1,\ldots,\gamma_N)$ such that
$\alpha=\gamma_1\v a_1+\cdots+\gamma_N\v a_N$ and take as starting point is the formal Laurent series
$$\Phi_{L,\gamma}(v_1,\ldots,v_N)=\sum_{\v l\in L}{\v v^{\v l+\gamma}\over\Gamma(\v l+\gamma+\v 1)}$$
where we use the short-hand notation
$${\v v^{\v l+\gamma}\over\Gamma(\v l+\gamma+\v 1)}={v_1^{l_1+\gamma_1}\cdots v_N^{l_N+\gamma_N}
\over \Gamma(l_1+\gamma_1+1)\cdots\Gamma(l_N+\gamma_N+1)}.$$
Note that there is a freedom of choice in $\gamma$ by shifts over $L\otimes\bbbr$.
A priori this series is formal, i.e. there is no convergence. However by making proper choices for
$\gamma$ we do end up with series that have an open domain of convergence in $\bbbc^N$.

Choose a subset ${\cal I}\subset\{1,2,\ldots,N\}$ with $|{\cal I}|=N-r$ such that $\v a_i$
with $i\not\in {\cal I}$ are linearly independent. In \cite[Prop 1]{GKZ2} we find the following
proposition (albeit in a different formulation).
\begin{proposition}\label{simplexvol}
Define $\pi_{\cal I}:L\to \bbbz^{N-r}$ by $\v l\mapsto (l_i)_{i\in\cal I}$.
Then $\pi_{\cal I}$
is injective and its image is a sublattice of $\bbbz^{N-r}$ of index
$|\det(\v a_i)_{i\not\in {\cal I}}|$.
\end{proposition}

We denote $\Delta_{\cal I}=|\det(\v a_i)_{i\not\in {\cal I}}|$.
Choose $\gamma$ such that $\gamma_i\in\bbbz$ for $i\in {\cal I}$.
The formal solution series
$$\Phi=\sum_{\v l\in L}\prod_{i\in\cal I}{v_i^{l_i+\gamma_i}\over\Gamma(l_i+\gamma_i+1)}
\prod_{i\not\in\cal I}{v_i^{l_i+\gamma_i}\over\Gamma(l_i+\gamma_i+1)}$$
is now a powerseries because the summation runs over the polytope
$l_i+\gamma_i\ge0$ for $i\in\cal I$ and the other $l_j$ are dependent on $l_i,i\in{\cal I}$.
Terms where $l_i+\gamma_i< 0$ do not occur because $1/\Gamma(l_i+\gamma_i+1)$ is zero
when $l_i+\gamma_i$ is a negative integer. By slight abuse of language will
call the corresponding simplicial cone $l_i\ge0$ for $i\in\cal I$
the {\it sector of summation} with index $\cal I$.

Denote the resulting series expansion by $\Phi_{{\cal I},\gamma}$.
The following statement, which is a direct consequence of estimates using Stirling's formula for
$\Gamma$, says that there is a non-trivial region of convergence.
\begin{proposition}
Let $(\rho_1,\ldots,\rho_N)\in\bbbr^N$ be such that $\rho_1l_1+\cdots+\rho_Nl_N>0$
for all $\v l\in L$ with $\forall i\in{\cal I}:l_i\ge0$. Then
$\Phi_{{\cal I},\gamma}$ converges for all $\v v\in\bbbc^N$ with
$|v_i|=t^{\rho_i}$ for sufficiently small $t\in\bbbr_{>0}$.
\end{proposition}
A proof can be found for example in \cite{stienstra1}. An $N$-tuple $\rho$
such that  $\rho_1l_1+\cdots+\rho_Nl_N>0$
for all $\v l\in L$ with $\forall i\in{\cal I}:l_i\ge0$ will be called
a {\it convergence direction}.

The following statement is a direct Corollary of Proposition \ref{simplexvol}.
\begin{corollary}
With notations as above, the number of distinct choices modulo $L$ for $\gamma$ such that
$\forall i\in{\cal I}:\gamma_i\in\bbbz$ is $\Delta_{\cal I}$.
\end{corollary}

There is one important assumption we need in order to make this approach work. Namely
the garantee that not too many of the arguments $l_i+\gamma_i$ are a negative integer. Otherwise
we might even end up with a power series which is identically zero.
The best way to do is to impose the condition $\gamma_i\not\in\bbbz$ for $i\not\in{\cal I}$.
Geometrically, since $\alpha=\sum_{i=1}^N\gamma_i\v a_i\is\sum_{i\not\in {\cal I}}\gamma_i\v a_i\mod{\bbbz^r}$,
this condition comes down to the
requirement that $\alpha+\bbbz^r$ does not contain points in a face of the simplicial cone spanned by
$\v a_i$ with $i\not\in{\cal I}$. Unfortunately this is stronger than the requirement of non-resonance of $H_A(\alpha)$,
as faces of the individual simplicial cones, not necessarily on the boundary of $C(A)$, are involved. However, the condition of
non-resonance does turn out to be useful.

\begin{proposition}\label{fullsupport}
Let ${\cal I}$ be as above and suppose the system $H_A(\alpha)$ is non-resonant.
Then there exists an open cone $C$ in $L\otimes\bbbr$ such
the series $\Phi_{{\cal I},\gamma}$ has non-zero terms for all $\v l\in C$.
\end{proposition}

{\bf Proof}.
We will use the following observation. The $i$-th coordinate of $\v l\in L$ can be
considered as a linear form on $L$. We shall do so in this proof. Suppose we have a relation
$\sum_{i=1}^N\lambda_i l_i=0$ with $\lambda_i\in\bbbr$. Then there exists a linear
form $m$ on $\bbbr^r$ such that $m(\v a_i)=\lambda_i$ for $i=1,\ldots,N$.

Denote the set of indices $i$ for which $\gamma_i\not\in\bbbz$ by $R$. When $|R|=r$
all terms of $\Phi_{{\cal I},\gamma}$ are non-zero and our statement is proven. Suppose $|R|<r$. Then there exist
linear relations between the forms $l_i$ with $\lambda_i=0$ when $i\in R$. Consider the convex
hull $D$ of the forms $l_i$ for $i\not\in R$. Suppose this hull contains the trivial form $\v 0$.
In other words, there exists a relation with coefficients $\lambda_i\in\bbbr_{\ge0}$, not all zero, with $\lambda_i=0$
for all $i\in R$. Hence, by our observation, there exists a non-trivial form $m$ on $\bbbr^r$ such that
$m(\v a_i)=\lambda_i$ for all $i$. Hence we have found a non-trivial form with $m(\v a_i)\ge0$ for all
$i$ and $m(\v a_i)=0$ for $i\in R$. Therefore the $\bbbr_{\ge0}$-span of $\v a_i, i\in R$ is contained
in a face $F$ of $C(A)$.
Furthermore, $\alpha=\sum_{i=1}^N\gamma_i\v a_i\is\sum_{i\in R}\gamma_i\v a_i\mod{\bbbz^r}$. Hence
modulo $\bbbz^r$ the vector $\alpha$ lies in the face $F$. This contradicts our non-resonance assumption
and therefore the convex hull $D$ does not contain $\v 0$. Consequently, the set of inequalities
$l_i\ge0,\ i\not\in R$ has a polyhedral cone with non-empty interior as solution space in
$\bbbr^{N-r}$. The terms in $\Phi_{{\cal I},\gamma}$ with indices inside this cone are non-zero.

\qed
\medskip

The following Theorem was one of the discoveries made by Gel'fand, Kapranov and Zelevinsky.

\begin{theorem} Let $\rho$ be a convergence direction.
Then there exists a regular triangulation $T$ of $A$ such that the summation
sectors for which $\rho$ is a convergence direction are given by $J^c$ where
$J$ runs through the $(r-1)$-simplices in $T$.
\end{theorem}

In order to proceed it is now important that different choices of summation sectors
give independent series solutions. For this we require the following condition.

\begin{definition}\label{Tnonresonance}
For any subset $J\subset\{1,2,\ldots,N\}$ denote $A_J=\{\v a_j|j\in J\}$ and
let $Q(A_J)$ be the convex hull of the points in $A_J$.

Let $T$ be a regular triangulation of $A$. The parameter $\alpha$ will be called
$T$-nonresonant if $\alpha+\bbbz^r$ does not contain a point on the boundary of
any cone over a $(r-1)$-simplex $Q(A_J)$ with $J\in T$. We call the system $T$-resonant otherwise.
\end{definition}

Notice that the $T$-nonresonance condition implies the nonresonance condition.
Let us assume that $\alpha$ is $T$-nonresonant.
For any ${\cal I}=J^c$ with $J\in T$ and one of the  ${\rm Vol}(Q(A_J))$ choices of $\gamma$
we get the series $\Phi_{{\cal I},\gamma}$.

\begin{theorem}
Under the $T$-nonresonance condition the power series solutions just constructed
form a basis of solutions of $H_A(\alpha)$.
\end{theorem}

{\bf Proof}.
To show that the solutions are independent it suffices to show that for any two
distinct summation sectors ${\cal I}$ and ${\cal I}'$ the
values of $\gamma_1,\ldots,\gamma_N$, as chosen in $\Phi_{\cal I}$ and $\Phi_{{\cal I}'}$,
are distinct modulo the lattice $L$. Suppose they are
not distinct modulo $L$. Then there exists an index $i\in{\cal I}'$, but $i\not\in{\cal I}$
such that $\gamma_i\in\bbbz$. But this is contradicted by our $T$-nonresonance assumption.

For every $J\in T$ we get ${\rm Vol}(Q(A_J))$ solutions by the different choices of
$\gamma$. Summing over $J\in T$ shows that we obtain
$\sum_{J\in T}{\rm Vol}(Q(A_J))={\rm Vol}(Q(A))$ independent solutions.

\qed

Given a regular triangulation we can consider the union of all summation domains in $L$.
More precisely, define $\supp(T)$ to be the convex closure of
$\cup_{J\in T}\{\v l\in L| l_i\ge0\ {\rm for\ all}\ i\in J^c\}$.
Then $\supp(T)$ will be the common support of all series $\Phi_{\cal I}$ with $I^c\in T$.
More precisely, denote the set of powerseries in $\v v$ with support in $\supp(T)$ by $\bbbc[[\v v]]_T$.
Note that this set forms a ring by the obvious multiplication. The coefficient ring $\bbbc$ can
be extended to the ring of finite linear combinations of powers $\v v^{\gamma}$ to get
the ring denoted by $\bbbc[\v v^{\gamma}][[\v v]]_T$. Note that the series constructed above
all belong to this ring. In the next section we further extend our coefficient ring to include
polynomials in $\log(v_i)$. This larger ring $\bbbc[\log(\v v),\v v^{\gamma}][[\v v]]_T$ is called
a Nilssen ring in \cite{SST}.

\section{T-resonant solutions}
In this section we assume that the system is $H_A(\alpha)$ is non-resonant, but
not necessarily $T$-nonresonant. In such a case
it is possible to write down a basis of solutions in $\bbbc[\log(\v v),\v v^{\gamma}][[\v v]]_T$.
This is done for example in \cite[Ch 3]{SST}.
We like to reproduce the proof from \cite{SST}, but in a slightly modified language.

Let $T$ be a regular triangulation of $Q(A)$. This time we assume the system to be
$T$-resonant when we specialise $\gamma$ to $\gamma^o$, say. Let
$$B_{\gamma^o}=\{J\in T| \gamma_i^o\in\bbbz\ {\rm for\ all\ }i\in J^c\}.$$
We say that the simplices $Q(A_J)$ with $J\in B_{\gamma^o}$ are resonating or {\rm in resonance}
with respect to $\gamma^o$. In case of $T$-nonresonance we would get $|B_{\gamma^o}|$ independent
series from the specialisation of $\gamma$ corresponding to $J\in B_{\gamma^o}$. Now we get
only one. So we have to find $|B_{\gamma^o}|-1$ additional series solutions. Just as in the
one variable case this will require the use of logarithms of the variables $v_i$

Let us denote $b=|B_{\gamma^o}|$. Choose $\alpha'$ such that $H_A(\alpha+\epsilon\alpha')$
is T-nonresonant for every sufficiently small $\epsilon\ne0$. The $b$ summation sectors $J^c$ with
$J\in B_{\gamma^o}$ now give rise to $b$ distinct specialisations
of the form $\gamma^o+\epsilon\gamma^{(i)}$ for $i=1,2,\ldots,b$ producing $b$ independent
solutions $\Phi_{\gamma^o+\epsilon\gamma^{(i)}}(\v v)$ of $H_A(\alpha+\epsilon\alpha')$. Multiply
each of these series by $\Gamma(\gamma^o+\epsilon\gamma^{(i)}+\v 1)$ to obtain the solutions
$$\Psi_i(\epsilon,\v v)=\sum_{\v l\in L}{\Gamma(\gamma^o+\epsilon\gamma^{(i)}+\v 1)
\over \Gamma(\v l +\gamma^o+\epsilon\gamma^{(i)}+\v 1)}\ \v v^{\v l+\gamma^o+\epsilon\gamma^{(i)}}.$$
Note that the coefficients are rational functions of $\epsilon$. Now expand
$$\v v^{\v l+\gamma^o+\epsilon\gamma^{(i)}}=
\sum_{n\ge 0}{\epsilon^n\over n!}(\gamma_1^{(i)}\log v_1+\cdots+\gamma_N^{(i)}\log v_N)^n.$$
Also expand the rational function
$\Gamma(\gamma^o+\epsilon\gamma^{(i)}+\v 1)/\Gamma(\v l+\gamma^o+\epsilon\gamma^{(i)}+\v 1)$
into a power series in $\epsilon$. We get
$$\Psi_i(\epsilon,\v v)=\sum_{n\ge0}{\epsilon\over n!}\Psi_i^{(n)}(0,\v v)$$
where $\Psi_i^{(n)}(\epsilon,\v v)$ denotes the $n$-th derivative of $\Psi_i(\epsilon,\v v)$
with respect to $\epsilon$. In particular, $\Psi_i(0,\v v)=\Gamma(\gamma^o+\v 1)\Phi_{\gamma^o}(\v v)$,
i.e. all $\epsilon$-series expansions $\Psi_i(\epsilon,\v v)$ have the same initial term.

Let $V_0$ be the $\bbbc$-vector space generated by the $\Psi_i(\epsilon,\v v)$. Its dimension is
$b$. There is a filtration $V_0\supset V_1\supset V_2\supset\cdots$ on $V_0$ defined by
$f(\epsilon,\v v)\in V_m$ if $f$ is divisible by $\epsilon^m$. Clearly $\dim(V_M)=0$ for
sufficiently large $M$. Let $f(\epsilon,\v v)\in V_m$. Then
$g(\v v)=\lim_{\epsilon\to0}\epsilon^{-m}f(\epsilon,\v v)$ is a solution of $H_A(\alpha)$.
This is clear for the box-operators $\Box_{\v l}$ since they are independent of $\epsilon$.
Let $Z_i$ be a homogeneity operator. Then $(Z_i-\alpha_i-\epsilon\alpha'_i)f(\epsilon,\v v)=0$.
Divide by $\epsilon^m$ and let $\epsilon\to0$. Then $(Z_i-\alpha_i)g(\v v)=0$, as desired.
Note that $g(\v v)\in\bbbc[\log(\v v),\v v^{\gamma}][[\v v]]_T$.

Let $b_j=\dim(V_j)$ for all $j$, in particular $b_0=b$. We choose a basis of $V_0$ as follows. Take
$b_0-b_1$ elements $f_{b_0},\ldots,f_{b_1+1}$ of $V_0$ which are linearly independent modulo $V_1$.
Choose $b_1-b_2$ elements $f_{b_1},\ldots,f_{b_2+1}\in V_1$ which are independent modulo $V_2$, etc.
We say that $f_i$ has weight $w$ if $f\in V_w$ and $f\not\in V_{w+1}$.
Divide $f_i$ by $\epsilon^w$ and let $\epsilon\to0$. Denote the limit by $g_i(\v v)$.
By construction elements $g_i(\v v)$ coming from $f_i$ of the same weight are linearly independent.
Elements $g_i(\v v)$ coming from $f_i$ with distinct weights are independent
because the series expansion have different degrees in the $\log(v_i)$. Hence the series $g_i(\v v)$
provide the desired $b$ independent solutions of $H_A(\alpha)$. Thus we obtained the Theorem
of Saito-Sturmfels-Takayama \cite[Thm 3.5.1]{SST} for the case of non-resonant systems (in their book
the author also produce bases of resonant systems).

\begin{theorem}[Saito-Sturmfels-Takayama]\label{sst}
Suppose $H_A(\alpha)$ is non-resonant. For any regular triangulation of $Q(A)$ there exists a
space of solutions to $H_A(\alpha)$ in the ring $\bbbc[\log(\v v),\v v^{\gamma}][[\v v]]_T$ of
$\bbbc$-dimension ${\rm Vol}(A)$.
\end{theorem}

By a Theorem of Adolphson \cite[Corollary 5.20]{adolphson} the rank of $H_A(\alpha)$ equals ${\rm Vol}(A)$
when the system is non-resonant. Hence we get the following.

\begin{corollary}\label{basis}
When $H_A(\alpha)$ is non-resonant the system of solutions in Theorem \ref{sst} provides a basis of solutions to
$H_A(\alpha)$ in $\bbbc[\log(\v v),\v v^{\gamma}][[\v v]]_T$.
\end{corollary}

\section{Non-resonant systems}\label{proofn2i}
In this section we prove Theorem \ref{nonresonant2irreducible}. Suppose we have a non-resonant system
and an operator $P\in K[\D]$ which annihilates a non-trivial solution $f$ in the solution space of $H_A(\alpha)$.

First we show the existence of such an $f$ which is of the form a power series of the type $\Phi_{\gamma}$,
as in the previous two sections. Fix a convergence direction $\rho_1,\ldots,\rho_N$ and let $T$ be
the corresponding regular triangulation of $Q(A)$.

Corollary \ref{basis} provides
a basis of solutions in $\bbbc[\log(\v v),\v v^{\gamma}][[\v v]]_T$. Consider these solutions as analytic functions
on an open neighbourhood of the set $V$ given by $|v_1|=t^{\rho_1},\ldots,|v_N|=t^{\rho_N}$ for $t$ sufficiently small.
The fundamental group
$\pi_1(V)$ is generated by $v_j=t^{\rho_j}e^{2\pi ix},\ x\in[0,1]$ for any $j$ and $v_i$ fixed for all $i\ne j$.
The corresponding monodromy group is an abelian group and so is its restriction to the common solution space of $H_A(\alpha)$
and $P(f)=0$. Since the monodromy group is abelian, there exists a one-dimensional invariant subspace. The character,
with which $\pi_1(V)$ acts on this space, uniquely determines a solution of the form $\Phi_{\gamma}$.

In the terminology of \cite[Thm 2.7]{mat} the solution $\Phi_{\gamma}$ is a fully supported
solution by virtue of Proposition \ref{fullsupport}. Theorem 2.7 of \cite{mat} implies that
the operator $P$ lies in ${\cal H}_A(\alpha)$. Hence we conclude that $H_A(\alpha)$ is irreducible.
\qed

\end{document}